\documentclass{amsart}

\newtheorem{Theorem}{Theorem}
\newtheorem{Hypothesis}[Theorem]{Hypothesis}
\newtheorem{Definition}[Theorem]{Definition}
\newtheorem{Remark}[Theorem]{Remark}
\newtheorem{Proposition}[Theorem]{Proposition}
\newtheorem{Corollary}[Theorem]{Corollary}

\begin{document}

\title[SOC via SPDEs]{Self-organized criticality via stochastic partial differential equations}

\author{Viorel Barbu}
\address{Institute of Mathematics, ``Octav Mayer'', Iasi, Romania}

\author{Philippe Blanchard}
\address{Faculty of Physics, University of Bielefeld, Germany}

\author{Giuseppe Da Prato}
\address{Scuola Normale Superiore di Pisa, Italy}

\author{Michael R\"ockner}
\address{Faculty of Mathematics, University of Bielefeld, Germany, \\
  Departments of Mathematics and Statistics, Purdue University,  U. S. A.}

\maketitle

\begin{abstract}
Models of self-organized criticality which can be described 
as singular diffusions with or without 
(multiplicative) Wiener forcing
term (as e.g.\ the Bak/Tang/Wiesenfeld- and Zhang-models) are 
analyzed.
Existence and uniqueness of nonnegative strong solutions are 
proved.
Previously numerically predicted transition to the critical 
state in 1-D is confirmed  by a rigorous proof that this indeed
happens in finite time with high probability.
\end{abstract}

\section{Introduction}
Within the past twenty years the notion of self-organized criticality 
(SOC) has become a new paradigm for the explanation of a huge variety 
of phenomena in nature and social sciences. 
Its origin lies in the attempt to explain the widespread appearance 
of power-law-statistics for characteristic events. 
In this paradigm an external perturbation may induce a chain reaction 
or avalanche in the system. 
Furthermore, a stationary state, the SOC-state, is reached where the 
average incoming flux is balanced by the average outgoing flux. 
This phenomenon was quite unexpected since attaining the critical 
state of a thermodynamic system usually requires a fine tuning of 
some control parameter, which is absent in the definition of the 
SOC models. 
Several models have been proposed to mimic this mechanism including 
the sand pile BTW-model
\cite{bak-jang-wiesenfeld-1,bak-jang-wiesenfeld-2}
and the Zhang-model \cite{zhang}. 
The presence of thresholds in the definition of the dynamics implies 
that the energy can be accumulated locally, eventually generating a 
chain reaction which may transport energy on arbitrary large scales. 

The literature on SOC is vast. We refer e.g.\ to 
\cite{bak-jang-wiesenfeld-1,bak-jang-wiesenfeld-2}, 
\cite{bantay-janosi}, \cite{blanchard-cess-krue}, 
\cite{jensen-1988}, \cite{turcotte}, 
\cite{carlson-chay-1, carlson-chay-2}, 
\cite{hentschel-family}, \cite{janosi-kertesz}, 
\cite{diaz-guil}, \cite{cafiero-loret-piet}, 
\cite{gia-diaz}, \cite{hwa-kardar}, 
\cite{grin-lee-sach}, \cite{zhang} for various studies. 
In \cite{bantay-janosi} it was proposed to describe this 
phenomenon, e.g.\ in the case of the avalanche dynamics
in the BTW- (see \cite{bak-jang-wiesenfeld-1,bak-jang-wiesenfeld-2}) and Zhang- 
(see \cite{zhang}) models, by a singular diffusion.
In the absense of noise the density $\varrho(t,\xi)$, 
$t\geq 0$, $\xi\in\mathbb R^d$, of this diffusion is formally
described by the evolution equation
\begin{equation}\label{eq-1--1}
  \frac{\partial}{\partial t}\varrho(t,\xi)
  = \Delta\Psi(\varrho(t,\xi)) ,
\end{equation}
where $\Psi(\varrho) := f(\varrho) H(\varrho-\varrho_c)$,
$\varrho_c \geq 0$, $H$ is the Heaviside function and $f(\varrho)=\text{const.}$
in the BTW-model and $f(\varrho)=\varrho$ in the Zhang-model.
Setting 
$D(\varrho) = f'(\varrho)H(\varrho-\varrho_c)
  + f(\varrho)\delta(\varrho-\varrho_c)$, 
the equation turns into
\begin{equation}\label{eq-1-0}
\begin{split}
  &\frac{\partial}{\partial t}\varrho(t,\xi)
  = \nabla\cdot [D(\varrho(t,\xi))\,\nabla\varrho(t,\xi)] \\
  &\quad= D'(\varrho(t,\xi))\, |\nabla\varrho(t,\xi)|^2
    + D(\varrho(t,\xi))\,\Delta\varrho(t,\xi).
\end{split}
\end{equation}
To discuss the problem in a heuristic way, let us consider a 
smooth version of $H$, for example
\begin{equation*}
  H_\varepsilon(\varrho) 
  = \frac 12 + \frac 1\pi \arctan (\tfrac \varrho\varepsilon)
\end{equation*}
(or the mathematical more convenient one in \eqref{e2.3a} 
below) 
and the corresponding function 
$D_\varepsilon = D_\varepsilon(\varrho)$.
Since $H_\varepsilon(\varrho-\varrho_c)$ is convex left of 
$\varrho_c$ and concave right of $\varrho_c$, we see that e.g.\ in 
the BTW-model $D'_\varepsilon(\varrho)|\nabla\varrho|^2$ is
positive if $\varrho$ is below $\varrho_c$ and negative if
$\varrho$ is above $\varrho_c$, i.e.\ according to
\eqref{eq-1-0}, $\varrho$ is ``pushed'' towards the critical
value $\varrho_c$.
This has been predicted numerically in 1-D by Bantay and Janosi in \cite{bantay-janosi}.

In \cite{diaz-guil} (see also \cite{gia-diaz}) Diaz--Guilera pointed out that
it is more realistic to consider \eqref{eq-1--1} perturbed by (an additive) noise to
model a random amount of energy put into the system varying all over the underlying
domain. So, \eqref{eq-1--1} turns into a stochastic partial differential equation (SPDE).
In \cite{carlson-swindle} (see also the references therein) based on numerical tests,
Carlson and Swindle observed that in the presence of such a noise the self-organized
behaviour does not necesssarily occur, i.e. the system fails to converge to the
critical value $\varrho _c$.

The purpose of this note is to provide rigorous mathematical proofs for the above
phenomena. First, we sketch the proofs for existence and uniqueness of solutions to
\eqref{eq-1--1} perturbed by noise (more precisely for multiplicative noise, so
that positivity of initial data is preserved). Second, we prove that at least
in 1-D we have convergence to $\varrho_c$ in finite time in the deterministic
case (confirming the numerical results in \cite{bantay-janosi}) and convergence to $\varrho_c$
with high probability in the stochastic case. In regard to \cite{carlson-swindle} one
can probably not achieve more, but so far we failed to prove that this probability is
really not equal to $1$.

Let us introduce our framework, where we switch to common notation
in SPDE, i.e.\ replace $\varrho(t,\xi)$ by $X(t,x)(\xi)$
with $x$ being the density at $t=0$, $t\geq 0$, $\xi\in\mathcal O$,
where $\mathcal O$ is
an open bounded domain in $\mathbb R^d, d=1,2,3,$ with smooth boundary  $\partial\mathcal O.$
The appropriate class of SPDE is then of the form
\begin{equation}
\label{e1.1}
\left\{\begin{array}{l}
dX(t)-\Delta\Psi(X(t))dt\ni \sigma(X(t))dW(t),\qquad\\
\hfill\mbox{\rm in}\;(0,\infty)\times  \mathcal O,\\
\Psi(X(t))\ni 0,\quad\mbox{\rm on}\; (0,\infty)\times \partial \mathcal O, \\
X(0,x)=x\quad\mbox{\rm on}\; \mathcal O,
\end{array}\right.
\end{equation}
 where $x$ is an initial datum, $\Psi:\mathbb R\to 2^\mathbb R$ a maximal monotone graph 
 and
$$
\sigma(X)dW=\sum_{k=1}^N \mu_kX d\beta_k\;e_k,\quad t\ge 0,
$$
is a random forcing term,
where   $\{e_k\} \subset L^2 (\mathcal O)$ is the eigen basis of the Laplacian $-\Delta $ on $\mathcal O$ with Dirichlet boundary conditions, 
$N\in \mathbb N \cup \{ + \infty\}$, $\mu_k$ are positive numbers and    $\beta_k$  independent   standard Brownian
motions on a filtered probability space $(\Omega,\mathcal F,\{\mathcal F_t\}_{t\ge 0},\mathbb P).$ 

Throughout this note we make the following assumptions:
\begin{Hypothesis} 
\label{h1.1}
\begin{enumerate}
\item[(i)] 
$\Psi(r)=\rho\;\mbox{\rm sign}\;r+\widetilde{\Psi}(r)$,  for   $r\in \mathbb R,$
where $0\in \Psi(0)$, $\rho>0$, $\widetilde{\Psi}:\mathbb R\to\mathbb R$ is  Lipschitzian, 
$\widetilde{\Psi}\in C^1(\mathbb R\setminus\{0\})$  and  for some $\delta>0$ it  
satisfies $\widetilde{\Psi}'(r)\ge \delta $ for all $r\in \mathbb R\setminus\{0\}$.
\item[(ii)] If $N=\infty$, the sequence   $\{\mu_k\}$  is such that
$$
\sum_{k=1}^\infty\mu^2_k\lambda^2_k<+\infty,
$$
where $\lambda_k$ are the eigenvalues of $-\Delta$.
\end{enumerate}
\end{Hypothesis} 
A typical example is given by 
\[\Psi(r)= \psi _0 (r) + c , \; r \in\mathbb R,\]
where 
\[\psi_0(r):= \begin{cases}
               \alpha_1 r ,&\quad r> 0\\ 
               [-\varrho, \varrho], &\quad r =0 \\
               \alpha_2 r , & \quad r <0
              \end{cases}
\]
and $\alpha_1, \alpha_2 > 0, \, \varrho \geq 0, \, c\in [-\varrho , \varrho]$ are constants.
\bigskip

The following notations will be used.
$L^p(\mathcal O),\;p\ge 1,$ is the usual space of $p$-integrable functions on $\mathcal O$ with norm $|\cdot|_p$. The scalar product in $L^2(\mathcal O)$ and the duality induced by the pivot space $L^2(\mathcal O)$
will be denoted by $\langle \cdot, \cdot   \rangle_2$.
$H^1_0(\mathcal O)\subset L^2(\mathcal O)$ is the first order Sobolev space on $\mathcal O$ with
zero trace on the boundary. 
For a fixed measure space $(E , \mathcal E, m)$, a Banach space $B$ and $p\in [1,\infty]$ 
we denote the space of all with respect to the measure $m$ $p$-integrable maps from $E$ to $B$ by 
$L^p(E;B)$. 

In the following by $H$ we shall denote  the distribution space
$$
H=H^{-1}(\mathcal O)=(H^{1}_0(\mathcal O))'
$$
endowed with the scalar product and norm defined by
$$
\langle u,v   \rangle=\int_\mathcal O (-\Delta)^{-1}u(\xi)v(\xi)d\xi,\quad |u|_{-1}=\langle u,u   \rangle^{1/2}.
$$
We recall that the operator $x\to -\Delta \Psi(x)$ with the domain
$$
\{x\in L^1( \mathcal O)\cap H: \exists \eta\in H^1_0(\mathcal O) \text{ s.th.\ $\eta\in\Psi(x)$ a.e.\ on $\mathcal O$} \},
$$
is maximal monotone in $H$ (see e.g. \cite{barbu-1976}) and so the distribution space $H$  is the natural functional setting for equation \eqref{e1.1}. However, the general existence theory of infinite dimensional stochastic equations in Hilbert space with nonlinear maximal monotone operators (see e.g. \cite{daprato-zabczyk-96}, \cite{prevot-roeckner-07}) is not directly applicable and so a direct approach must be used.

\section{Existence, uniqueness and positivity}
\begin{Definition}
\label{d2.1}
Let $x\in H$. An $H$-valued   continuous $\mathcal F_t$-adapted process $X=X(t,x)$ is called a solution to \eqref{e1.1} 
 on $[0,T]$ if for some $p\in [1,\infty[$
$$
X\in L^{p} (\Omega\times(0,T)\times \mathcal O)\cap L^2(0,T;L^2(\Omega,H)),
$$
 and there exists $\eta\in 
L^{p} (\Omega\times(0,T)\times \mathcal O)$ such that $\mathbb P$-a.s.
\begin{equation}
\label{e2.1}
\begin{array}{lll}
&\langle X(t,x),e_j\rangle_2 &
\\=&\langle x,e_j\rangle_2 +\int_0^t\int_\mathcal O\eta(s,\xi)\Delta e_j(\xi) d\xi ds\\
 &\quad +\sum_{k=1}^\infty\mu_k\int_0^t\langle X(s,x)e_k,e_j\rangle_2 d\beta_k(s),\quad \\
 &\hfill \forall\;j\in \mathbb N,\;t\in [0,T],
\end{array}
\end{equation}
\begin{equation}
\label{e2.2}
\eta\in \Psi(X)\;\quad\mbox{\rm a.e. in }\;\Omega\times (0,T)\times \mathcal O.
\end{equation}
\end{Definition}
Below for simplicity we often write $X(t)$ instead of $X(t,x)$.

\begin{Theorem}
\label{t2.2}
Under Hypothesis $\ref{h1.1}$ for each $x\in L^p( \mathcal O)$, $p\geq 4$  
there is a unique solution $X$
to \eqref{e1.1}. Moreover, if $x$ is nonnegative a.e. in $\mathcal O$ then $\mathbb P$-a.s.
$$
X(t,x)(\xi)\ge  0, \quad\mbox{\it for a.e.}\;(t,\xi)\in (0,\infty)\times \mathcal O.
$$
\end{Theorem}

{\noindent\bf Sketch of Proof}. 
Consider the approximating equation
\begin{equation}
\label{e3.1}
\left\{\begin{array}{l}
dX_\lambda(t)-\Delta(\Psi_\lambda(X_\lambda(t))
dt=\sigma(X_\lambda(t))dW(t),\ \\
X_\lambda(0,x)=x,
\end{array}\right.
\end{equation}
where $\lambda>0$,
\begin{equation}\label{e2.3a}
\Psi_\lambda(r):=\rho\;(\mbox{\rm sign})_\lambda(r)+\widetilde{\Psi}(r),\quad r\in \mathbb R,
\end{equation} 
 $$
 (\mbox{\rm sign})_\lambda(r):=
\left\{\begin{array}{l}
\protect
1\quad\mbox{\rm if}\;  r>\lambda\\
\frac{r}\lambda\quad\mbox{\rm if}\;  r\in[-\lambda,\lambda]\\
-1\quad\mbox{\rm if}\;  r<-\lambda.
\end{array}\right.
$$

By \cite[Theorem 2.2]{barbu-dap-roe} (applied with $m=1$) equation \eqref{e3.1} has a unique solution
$$
X_\lambda\in L^{2}(\Omega\times (0,T)\times\mathcal O)\cap L^2(\Omega,C([0,T];H))
$$
in the sense of Definition \ref{d2.1} which is nonnegative, if so is $x$. 
Here as usual the space of continuous $H$-valued paths
 $C([0,T];H)$ is equipped with the supremum norm.

By Ito's formula for $\alpha>0$ large enough it follows that for all $\lambda,\mu\in (0,1)$ and $t\in[0,T]$
\begin{equation}
\label{e3.11}
\begin{array}{l}
\frac12\;|X_\lambda(t)-X_\mu(t))|_{-1}^2 e ^{-\alpha t}\\
\le C \max\{\lambda,\mu\}\int_0^t\int_\mathcal O
\Big(|\Psi_\lambda(X_\lambda(s))|^2+|X_\lambda(s)|^2\\
+|\Psi_\mu(X_\mu(s))|^2 +|X_\mu(s)|^2\Big)e^{-\alpha s}d\xi\;ds\\
 +\int^t_0 e^{-\alpha s}\langle X_\lambda (s)
- X_\mu (s),
\sigma(X_\lambda (s) - X_\mu (s))dW(s)\rangle_2 .
\end{array} 
\end{equation}
Hence by the Burkholder-Davis-Gundy inequality 
$\{X_\lambda \}$ is a Cauchy net in $L^2(\Omega; C([0,T], H))$ and by a standard technique from stochatic PDE one 
shows that the limit $X$ is the desired solution to \eqref{e1.1} (cf. \cite{barbu-dap-roe-2}).
\hfill$\Box$ 

\begin{Remark}
 One can also show (see \cite[Prop. 3.4]{barbu-dap-roe-2}) that 
 $X, \, X_\lambda \in L^2 (0,T; L^2 (\Omega, H_0^1 (\mathcal O)))$, that
 \[\lim_{\lambda \to 0} \mathbb E \int _0^T | X_\lambda - X| ^2 _{L^2 (\mathcal O)} dt =0,\] 
 where $\mathbb E$ denotes expectation with respect to $\mathbb P$, and that both $X$ and $X_\lambda$ have 
 continuous paths in $L^2(\mathcal O)$. Theorem \ref{t2.2} is true for more general not linear growing $\Psi$
 (see \cite{barbu-dap-roe-2}).
\end{Remark}
\section{Extinction in finite time and self-organized criticality}
In this section we assume $N<\infty$. Let $\tau$ be the stopping time
$$ 
\tau=\inf\{t\ge 0:\;|X(t,x)|_{-1}=0\},
$$ 
 where $X(t,x), t\ge 0,$ is the solution from Theorem \ref{t2.2}.
 \begin{Proposition}
\label{l4.1}
$$
X(t,x)=0\quad\mbox{\rm for} \; t\ge \tau,\;\;\mathbb P\mbox{\rm -a.s.}.
$$
\end{Proposition}
{\noindent\bf Sketch of Proof}. 
For simplicity we consider the case with $\varrho \equiv 1$.
Define
$$
\mu(t):=-\sum_{k=1}^N\mu_ke_k\beta_k(t),\quad t\in [0,T],\quad
\tilde\mu:=\sum_{k=1}^N\mu^2_ke^2_k
$$
and
$$
Y(t):=e^{\mu(t)}X(t),\quad t\ge 0.
$$
Then by Ito's product rule $Y$ satisfies $\mathbb P$-a.s. the following ordinary PDE
\[\frac{dY(t)}{dt}= e^{\mu(t)} \Delta \eta(t) - \frac 12 \tilde \mu Y(t),\quad t \geq 0.\]
Setting $Y_\lambda := e^{\mu }X_\lambda$ we consider the approximating equation
\begin{equation}
\label{e4.4}
\frac{dY_\lambda(t)}{dt}=e^{\mu(t)}\Delta \eta_\lambda(t)-\frac12\;\tilde\mu(t)Y_\lambda(t),\quad \;t\ge 0,
\end{equation}
where
\begin{equation}\label{e3.1a}
\eta_\lambda(t)=\Psi_\lambda(X_\lambda(t)) \in H^1_0(\mathcal O).
\end{equation} 
Hence
\begin{equation}
\label{e4.7}
\begin{split}
&\left<\frac{dY_\lambda(t)}{dt},Y_\lambda(t)\right>_2\\
=&
\left<\eta_\lambda(t),\Delta(e^{\mu(t)}Y_\lambda(t))\right>_2 
  -\frac12 \langle \tilde\mu(t)Y_\lambda(t),Y_{\lambda}(t) \rangle_2,
\end{split}
\end{equation}
where by \eqref{e3.1a}, \eqref{e2.3a} and integrating by parts we have
$$
\begin{array}{l}
\displaystyle \langle\eta_\lambda(t),\Delta(e^{\mu(t)}Y_\lambda(t))   \rangle_2
\\
\displaystyle = -\frac1\lambda\;\int_{\{ e^{-\mu (t)}|Y_\lambda (t)| < \lambda\}}(|\nabla Y_\lambda(t)|^2\\
\qquad\qquad\qquad-|Y_\lambda(t)|^2 
\;|\nabla\mu(t)|^2) d\xi\\
\displaystyle  \quad - \int_\mathcal O\widetilde{\Psi}'(e^{-\mu(t)}Y_\lambda(t))\\
\qquad\qquad\qquad(|\nabla Y_\lambda(t)|^2-|Y_\lambda(t)|^2 
\;|\nabla\mu(t)|^2)  d\xi.
\end{array}
$$
This yields
\begin{equation}
\label{e4.7'}
\langle\eta_\lambda(t),\Delta(e^{\mu(t)}Y_\lambda(t))   \rangle_2
\le C\left(|Y_\lambda(t)|^2_2+\lambda\right).
\end{equation}
Hence \eqref{e4.7} and Gronwall's lemma imply
$$
|Y_\lambda(t)|^2_2\le e^{C(t-s)}\left(|Y_\lambda(s)|^2_2+C\lambda T\right),\quad t\geq s.
$$
Now letting $\lambda\to 0$  we get
\begin{equation}
\label{e4.8}
|Y(t)|^2_2\le e^{C(t-s)}|Y(s)|^2_2,\quad t\geq s.
\end{equation}
Taking in \eqref{e4.8} $s=\tau$ we get $Y(t)=0$ for all $t\ge \tau$ as claimed. 
\hfill $\Box$\vspace{3 mm}

For proving the extinction result we need $\mathcal O\subset \mathbb R$, i.e. $d=1$. To be more specific let $\mathcal O=(0,\pi)$. Then
$e_k(\xi)=\sqrt{\frac2\pi}\;\sin k\xi,\quad \xi\in [0,\pi]$, $\lambda_k=k^2$ and $L^1(0,\pi)\subset H$ continuously, so
\begin{equation}
\label{e4.11}
\gamma=\inf\left\{ \frac{|x|_{L^1}}{|x|_{-1}}:\;x\in L^1(0,\pi)  \right\}>0.
\end{equation}

\begin{Theorem}\label{t3.2}
  Consider the equation
\begin{equation}
\label{e4.19}
\left\{\begin{array}{l}
dX(t)-\Delta(\rho\;\mbox{\rm sign}\;(X(t)-x_c)+\widetilde{\Psi}(X(t)-x_c))dt\\
\hspace{20mm}\displaystyle\ni \sigma (X(t)-x_c)\sum_{k=1}^N\mu_ke_kd\beta_k,\quad t\ge 0, \\
\rho\;\mbox{\rm sign}\;(X(t)-x_c) + \widetilde \Psi (X(t) - x_c)\ni 0,\;\mbox{\rm on}\;\partial (0,\pi),\\
X(0,x)=x.
\end{array}\right.
\end{equation}
where $x_c\in \mathbb R$.

Assume that
$$
|x-x_c|_{-1}<\rho\gamma C_N^{-1},
$$
where $C_N :=\frac\pi{4}\;\sum_{k=1}^N(1+k)^2\mu_k^2$  and $\gamma$ is as in \eqref{e4.11}.
Then for all $t>0$
\begin{equation}\label{e4.20}
\mathbb P(\tau_c \leq t)
  \geq 1 -\frac{|x-x_c|_{-1}}{\rho \gamma\, \bigl(\int_0^t e^{-C_N s}ds\bigr)} ,
\end{equation}
where
\[
\begin{split}
\tau_c&=\inf\{t\ge 0:\;|X(t)-x_c|_{-1} =0  \}\\
&=\sup\{t\ge 0:\;
|X(t)-x_c|_{-1}>0\}
\end{split}
\]
and $X-x_c$  is the solution from Theorem \ref{t2.2}.
\end{Theorem}

{\noindent\bf Sketch of Proof}. 
For simplicity we assume $x_c=0$. Since $\langle X(t),\tilde \Psi(X(t))\rangle_{2}\geq 0$,
$0\in \Psi(0)$. An application of Ito's formula for $\varphi_\varepsilon (|X|^2_{-1})=(|X|_{-1}^2 + \varepsilon ^2)^{1/2}$ yields
$\mathbb P$-a.s.
\begin{equation}
\label{e4.17}
\begin{array}{l}
\displaystyle\varphi_\varepsilon(|X(t)|^2_{-1})+\gamma \rho\int_0^{\min (t\wedge\tau)}\frac{|X(s)|_{-1}}{(|X(s)|_{-1}^2+\varepsilon^2)^{1/2}}ds\\
\displaystyle\le\varphi_\varepsilon(|x|^2_{-1})+C_N
\int_0^t  \frac{|X(s)|^2_{-1}}{(|X(s)|_{-1}^2+\varepsilon^2)^{1/2}}ds\\
\displaystyle+2\int_0^{\min (t\wedge\tau)}\langle\sigma(X(s))dW(s),\varphi'_\varepsilon(|X(s)|^2_{-1})X(s) \rangle.
\end{array} 
\end{equation}
Now, letting $\varepsilon$ tend to zero we get $\mathcal P$-a.s.\ for $t\geq 0$
\begin{equation}
\label{e4.18}
\begin{array}{l}
\displaystyle |X(t)|_{-1}+\gamma\rho\, {\min (t\wedge\tau)}\\
\le|x|_{-1}+C_N \int_0^t  |X(s)|_{-1}ds\\
\displaystyle+2\int_0^t 1_{[0,\tau]} (s)\, \langle \sigma(X(s))dW(s),X(s)|X(s)|^{-1}_{-1}  \rangle.
\end{array} 
\end{equation}
Hence by a standard comparison result
\begin{align*}
 |X(t)| _{-1} + \gamma \rho \int _0^t e^{C_N(t-s)}1_{[0,\tau]}(s) ds \\
 \leq e^{C_N t}|x|_{-1} 
 + \int_0^t\langle \sigma (X(s)) dW(s), X(s) | X(s)|_{-1}^{-1} \rangle .
\end{align*}
Taking expectation we get
\[\int _0^t e^{-C_N s} \mathbb P (\tau > s ) ds \leq \frac{|x|_{-1}}{\gamma \rho}.\]
Writing $\mathbb P(\tau> s )=1- \mathbb P (\tau \leq s)$ we deduce \eqref{e4.20}. \hfill $\Box$ 

\begin{Corollary}
  \label{corollary-7}
  If in the situation of the above Theorem~\ref{t3.2}, the noise is zero, i.e.\
  $C_N=0$, then 
  \[
    \tau_c\leq\frac{|x-x_c|_{-1}}{\varrho\gamma}.
  \]
\end{Corollary}
\noindent\textbf{Proof.}
The assertion follows from
\eqref{e4.18}. \hfill$\Box$

\textbf{Acknowledgments:}
This work has been supported in part by
the CEEX Project 05 of Romanian Minister of Research,
the DFG-International Graduate School ``Stochastics and Real World Models'', 
the SFB-701, NSF-Grants 0603742, 0606615 as well as the BiBoS-Research Center,
the research programme ``Equazioni di Kolmogorov'' from the Italian
``Ministero della Ricerca Scientifica e Tecnologica''
and ``FCT, POCTI-219, FEDER''.

\end{document}